\documentclass[a4paper, BCOR=0.0mm, DIV=calc]{scrartcl}
\usepackage[T1]{fontenc}
\usepackage[utf8]{inputenc}
\usepackage[ngerman]{babel}
\usepackage[osf,sc]{mathpazo}
\usepackage{textcomp}
\usepackage{microtype}
\usepackage{amsmath, amssymb, amsfonts}
\KOMAoptions{twoside=false, twocolumn=false, headinclude=false, footinclude=false, mpinclude=false, pagesize=auto}
\recalctypearea
\newcommand{\ebinom}[2]{\left(\frac{#1}{#2} \right)}
\begin{document}
\boldmath
\title{Über den allgemeinen Term der hypergeometrischen Reihen\footnote{
Originaltitel: "`De termino generali serierum hypergeometricarum"', erstmals publiziert in "`\textit{Nova Acta Academiae Scientarum Imperialis Petropolitinae} 7, 1793, pp. 42-82"', Nachdruck in "`\textit{Opera Omnia}: Series 1, Volume 16, pp. 139 - 162 "', Eneström-Nummer E652, übersetzt von: Alexander Aycock, Textsatz: Artur Diener, im Rahmen des Projektes "`Eulerkreis Mainz"' }}
\unboldmath
\author{Leonhard Euler}
\date{}
\maketitle
\paragraph{§1}
Ich nenne hier in Übereinstimmung mit Wallis die Reihen hypergeometrisch, deren Terme nach  ununterbrochen multiplizierten Faktoren fortschreiten, während die Faktoren selbst eine arithmetische Progression beschreiben, von welcher Art die bekannteste Reihe $ 1,\, 2,\, 6,\, 24,\, 120,\, 720,\, \mathrm{etc.} $ ist, deren Term, der dem Index $n$ entspricht, $ 1 \cdot 2 \cdot 3 \cdot 4 \cdots n $ ist. Im Allgemeinen wird daher aus einer beliebigen arithmetischen Progression eine solche hypergeometrische Reihe
\[
	a,~ a(a+b),~ a(a+b)(a+2b),~ a(a+b)(a+2b)(a+3b),~ \mathrm{etc.}
\]
gebildet werden, deren Term, der dem Index $n$ entspricht, und der aus $n$ Faktoren besteht,
\[
	a(a+b)(a+2b)(a+3b) \cdots (a+(n-1)b)
\]
sein wird; daher, sooft $n$ eine ganze positive Zahl war, wird der Term, der selbigem entspricht, sehr leicht durch Faktoren solcher Art angegeben, sodass, wenn der Index $n$ die Abzisse einer beliebigen Kurve ausdrückt, ihre Ordinaten durch Terme der Reihe selbst ausgedrückt werden; es besteht natürlich kein Zweifel, dass auch den Abzissen, die entweder durch gebrochene Zahlen oder Wurzeln ausgedrückt wurden, bestimmten Ordinaten entsprechen, deren Größe sich daher keinesfalls bestimmen lässt; dafür wird ein Ausdruck solcher Art benötigt, der aus den Größen $a$, $b$ und $n$ gebildet wird, der immer den bestimmten Wert beschaffen soll, ob der Index $n$ eine ganze Zahl war oder eine gebrochene oder gar eine Wurzel.

\paragraph{§2}
Ich habe hypergeometrische Reihen solcher Art schon des Öfteren gründlich untersucht, wo ich hauptsächlich aus der Theorie der Interpolation der Wallis'schen Reihe $ 1,\, 2,\, 6,\, 24,\, 120,\, \mathrm{etc.} $ gefunden habe, dass der Term, der dem unbestimmten Index $n$ entspricht, so ausgedrückt werden kann, dass er gleich
\[
	\frac{1^{1-n} \cdot 2^n}{1+n} \cdot \frac{2^{1-n} \cdot 3^n}{2+n} \cdot \frac{3^{1-n} \cdot 4^n}{3+n} \cdot \frac{4^{1-n} \cdot 5^n}{4+n} \cdot \frac{5^{1-n} \cdot 6^n}{5+n} \cdot \mathrm{etc.}
\]
ist, welcher Ausdruck zwar ins Unendliche läuft, aber dennoch immer den bestimmten Wert ausdrückt, welcher Wert auch immer dem Index $n$ zugeteilt wird. Auf ähnliche Weise habe ich für die allgemeine oben erwähnte Reihe gezeigt, dass der allgemeine Term oder der dem unbestimmten Index $n$ entsprechende mit dem folgenden ins Unendliche laufende Produkt dargestellt werden kann:
\[
	\frac{a^n a^{1-n}(a+b)^n}{a+nb} \cdot \frac{(a+b)^{1-n}(a+2b)^n}{a+(n+1)b} \cdot \frac{(a+2b)^{1-n}(a+3b)^n}{a+(n+2)b} \cdot \mathrm{etc.}
\]
Die Methoden aber, die mich zu dieser Zeit zu diesen Formeln geführt haben, waren auf die Theorie der Interpolation beschränkt und vielleicht nicht so klar behandelt, dass sie hinreichend klar verstanden werden konnten; deshalb habe ich beschlossen, diese Erforschung aus der Gestalt dieser Reihen noch einmal herzuleiten und verständlich zu erklären.

\paragraph{§3}
Ich möchte also von der Wallis'schen Reihe aus anfangen, weil sich ja die Art der Schlussfolgerungen und Rechnungen in einem Spezialfall um Vieles besser erkennen lassen wird, als wenn ich diese sofort bei der allgemeinen Reihe anwenden wollte. Weil also der allgemeine Term, der dem Index $n$ entspricht, wie eine Funktion des Index' $n$ selbst betrachtet werden kann, möchte ich diesen auf hinreichend bekannte Weise durch $\Delta : n$ ausdrücken, wo $\Delta$ nicht eine Größe, sondern den Charakter der Funktion bezeichnet. Daher wird also, sooft $n$ eine positive ganze Zahl war,
\[
	\Delta : n = 1 \cdot 2 \cdot 3 \cdot 4 \cdots n
\]
sein, woraus eingesehen wird, dass für die folgenden Terme
\begin{align*}
	\Delta : (n+1) &= (n+1)\Delta : n \\
	\Delta : (n+2) &= (n+1)(n+2)\Delta : n \\
	\Delta : (n+3) &= (n+1)(n+2)(n+3)\Delta : n \\
	\mathrm{etc.}  &
\end{align*}
sein wird. Und weil ja in diesem ununterbrochenem Zugang der neuen Faktoren die Gestalt selbst der Reihe enthalten sein anzusehen ist, müssen diese letzen Formeln auch mit der Wahrheit verträglich sein, welche Werte auch immer dem Index $n$ zugeteilt werden. Weil $\Delta : \frac{1}{2}$ den Term bezeichnet, der dem Index $\frac{1}{2}$ entspricht und von dem bekannt ist, dass er durch die Quadratur des Kreises ausgedrückt wird, können aus diesem die folgenden so ausgedrückt werden:
\[
	\Delta : 1 + \frac{1}{2} = \frac{3}{2}\Delta : \frac{1}{2}, ~~ \Delta : 2\frac{1}{2} = \frac{3}{2} \cdot \frac{5}{2} \cdot \Delta : \frac{1}{2}, ~~ \Delta : 3\frac{1}{2} = \frac{3}{2} \cdot \frac{5}{2} \cdot \frac{7}{2} \cdot \Delta : \frac{1}{2}, ~~ \mathrm{etc.}
\]
Und auf diese Weise wird sich die Sache verhalten, welche andere Zahl auch immer für $n$ angenommen wird, wenn auch vielleicht sich der Wert auf keine Weise durch bekannte Maße ausdrücken lässt.

\paragraph{§4}
Nachdem also dieses Fundament gelegt wurde, auf welchem die ganze Natur dieser Reihen beruht, ist es notwendig ein anderes Prinzip festzusetzen, welches darin besteht, dass ins Unendliche fortgesetzte Reihen solcher Art schließlich mit einer geometrischen Progression verbunden werden, weil die neuen Faktoren, die weiter hinzukommen, als untereinander gleich betrachtet werden können. Wenn so $i$ eine unendliche Zahl bezeichnet und so $\Delta : i$ einen Term, der vom Anfang unendlich weit entfernt ist, bezeichnet, können die ihm folgenden Terme so beschafft werden:
\begin{align*}
\Delta : (i+1) &= (i+1)\Delta : i = i\Delta : i \\
\Delta : (i+2) &= (i+1)(i+2)\Delta : i = i\cdot i \Delta : i \\
\Delta : (i+3) &= (i+1)(i+2)(i+3)\Delta : i = i^3 \Delta : i \\
\mathrm{etc.}  &
\end{align*}
und so kann im Allgemeinen
\[
	\Delta : (i+n) = i^n \Delta : i
\]
gesetzt werden. Dies kann freilich auf den ersten Blick paradox erscheinen, weil ja die Formel $\Delta : i$ schon einen unendlich großen Wert hat; wo aber nur über das Verhältnis zwischen zwei oder mehreren Ausdrücken solcher Art gehandelt wird, wird sich sicherlich anstelle von $i+1$ und $i+2$ auch $i$ schreiben lassen. Ja, es wird sich sogar umgekehrt anstelle von $i$ auch $i+1$ oder $i+2$ oder im Allgemeinen $i+\alpha$ schreiben lassen, wobei $\alpha$ eine beliebige endliche Zahl bezeichnet, woher auch allgemein
\[
	\Delta : (i+n) = (i+\alpha)^n \Delta : i
\]
gesetzt werden.

\paragraph{§5}
Weil also im Allgemeinen
\[
	\Delta : i = 1 \cdot 2 \cdot 3 \cdot 4 \cdots i
\]
ist, wenn $i$ eine beliebige endliche Zahl bezeichnet, wenn auch vielleicht eine gebrochene, kann dieser Ausdruck nichtsdestotrotz als bekannt angesehen werden; dann aber wird, wie wir gesehen haben,
\[
	\Delta : (i+n) = (i+\alpha)^n\Delta : i
\]
sein. Auf dieselbe Weise wollen wir aber von der Formel $\Delta : n$ aus in gleicher Weise ins Unendliche vorgehen, und weil
\[
	\Delta : (n+1) = (n+1)\Delta : n
\]
ist und
\[
	\Delta : (n+2) = (n+1)(n+2)\Delta : n, ~ \mathrm{etc.}
\]
wird $ \Delta : (n+i) = (n+1)(n+2)(n+3) \cdots (n+i)\Delta : n$ sein; dort ist die Anzahl der Faktoren, mit welchen die Formel $ \Delta : n$ multipliziert wird, gleich $i$; aber tatsächlich wird der oben für $\Delta : (i+n)$ gegebene Ausdruck, wenn anstelle von $\Delta : i$ sein natürlicher Wert geschrieben wird,
\[
	\Delta : (i+n) = 1 \cdot 2 \cdot 3 \cdot 4 \cdots i(i+\alpha)^n
\]
sein, wo die Anzahl der Faktoren, mit denen die Formel $(i+\alpha)^n$ multipliziert wird, in gleicher Weise gleich $i$ ist.

\paragraph{§6}
Weil daher natürlich
\[
	\Delta : (n+i) = \Delta : (i+n)
\]
sein muss, muss, wenn die eine der gefundenen Formeln durch die andere geteilt wird und die Brüche, weil die Anzahl der Faktoren der beiden dieselbe ist, getrennt ausgedrückt werden, der Quotient natürlich der Einheit gleich hervorgehen und es wird natülich
\[
	1 = \frac{1}{n+1} \cdot \frac{2}{n+2} \cdot \frac{3}{n+3} \cdot \frac{4}{n+4} \cdots \frac{i}{n+i} \cdot \frac{(i+\alpha)^n}{\Delta : n}
\] 
sein. Aus dieser Gleichung lässt sich der wahre Wert der Formel $\Delta : n$ berechnen, der sogar immer gelten muss, ob der Index nun eine ganze Zahl ist oder anderer Art; wir werden natürlich
\[
	\Delta : n = \frac{1}{n+1} \cdot \frac{2}{n+2} \cdot \frac{3}{n+3} \cdot \frac{4}{n+4} \cdots \frac{i}{n+i} \cdots (i+\alpha)^n
\]
haben und so wie dieser Ausdruck unsere Aufgabe vollendet, wird es förderlich sein, in einigen Fällen, wo wir anstelle von $n$ ganze Zahlen, zumindest kleinere, angenommen haben, ihn gezeigt zu haben.
\begin{enumerate}
	\item Es sei daher $n=1$ und es wird daher 
	\[
		\Delta : 1 = \frac{1}{2} \cdot \frac{2}{3} \cdot \frac{3}{4} \cdot \frac{4}{5} \cdots \frac{i}{1+i} \cdot (i+\alpha)
	\]
	sein, wo nach Weglassen der sich aufhebenden Terme
	\[
		\Delta : 1 = \frac{1}{1+i} \cdot (i+\alpha)
	\]
	hevorgehen wird, wo natürlich wegen der unendlichen Zahl $i$ gilt, dass $\frac{i+\alpha}{1+i} =1$ ist, welcher Wert auch immer für $\alpha$ angenommen wird.
	\item Es sei $n=2$ und unser Ausdruck wird
	\[
		\Delta : 2 = \frac{1}{3} \cdot \frac{2}{4} \cdot \frac{3}{5} \cdot \frac{4}{6} \cdot \frac{5}{7} \cdot \frac{6}{8} \cdots \frac{i}{2+i}\cdot (i+\alpha)^2
	\]
	geben, wo nach Weglassen der sich streichenden Terme im Zähler nur die ersten zwei, im Nenner aber nur die letzten zwei übrig bleiben, so dass
	\[
		\Delta : 2 = \frac{1 \cdot 2 \cdot (i+\alpha)^2}{(1+i)(2+i)}
	\]
	hervorgeht, wo natürlich
	\[
		\frac{(i+\alpha)^2}{(1+i)(2+i)} = 1
	\]
	wird, so dass $\Delta : 2 = 1 \cdot 2$.
	\item Es sei $n=3$ und unser Ausdruck wird 
	\[
		\Delta : 3 = \frac{1}{4} \cdot \frac{2}{5} \cdot \frac{3}{6} \cdot \frac{4}{7} \cdot \frac{5}{8} \cdot \frac{6}{9} \cdot \frac{7}{10} \cdots \frac{i}{3+i}\cdot (i+\alpha)^3
	\]
	geben, wo im Zähler nur die ersten drei, im Nenner aber nun die letzten drei Faktoren übrigbleiben, so dass
	\[
		\Delta : 3 = \frac{1\cdot 2 \cdot 3 \cdot (i+\alpha)^3}{(1+i)(2+i)(3+i)}
	\]
	ist und daher
	\[
		\Delta : 3 = 1 \cdot 2 \cdot 3
	\]
	ist.
\end{enumerate}
Auf diese Art kann also die Gültigkeit dieses Ausdrucks für alle ganzen Zahlen gezeigt werden, und daher wird zugleich die Überlegung verstanden, warum jene beliebige Zahl $\alpha$ ohne Bedenken in die Rechnung eingeführt werden konnte, weil hier nur das Verhältnis zwischen zwei Unendlichkeiten in die Rechnung eingeht.

\paragraph{§7}
Wenn wir aber anstelle $n$ nicht ganze Zahlen annehmen würden, könnte aus dieser Form natürlich nichts für den Wert $\Delta : n$ gefolgert werden, weil im Zähler wie im Nenner unzählige Faktoren zurückblieben, unter welchen sogar unzählige selbst unendlich sein würden. Damit wir also dieser Unannehmlichkeit entgegen wirken, wollen wir, obwohl $i$ an sich eine unendliche Zahl bezeichnet, nichtsdestoweniger an ihrer Stelle nacheinander die natürlichen Zahlen $ 1,\, 2,\, 3,\, 4,\, \mathrm{etc.} $ schreiben und wir werden die folgenden Formeln erhalten
\begin{align*}
\mathrm{I.}~~~ \Delta : n &= \frac{1}{n+1}(1+\alpha)^n \\
\mathrm{II.}~~ \Delta : n &= \frac{1}{n+1} \cdot \frac{2}{n+2} (2 + \alpha)^n \\
\mathrm{III.}~ \Delta : n &= \frac{1}{n+1} \cdot \frac{2}{n+2} \cdot \frac{3}{n+3}(3+\alpha)^n \\
\mathrm{IV.}~~ \Delta : n &= \frac{1}{n+1} \cdot \frac{2}{n+2} \cdot \frac{3}{n+3} \cdot \frac{4}{n+4}(4 + \alpha)^n \\
\mathrm{etc.} & &
\end{align*} 
wo ersichtlich ist, dass diese Formeln sich immer mehr der Wahrheit nähern müssen, je weiter sie fortgesetzt werden, weil ja, nachdem sie ins Unendliche fortgesetzt wurden, man schließlich zum wahren Wert von $\Delta : n$ selbst gelangen muss.

\paragraph{§8}
Weil eine beliebige dieser Formeln die vorhergehende entweder ganz oder teilweise involviert, wollen wir eine beliebige durch ihre vorhergehende teilen; und wir werden
\begin{align*}
\frac{\mathrm{II}}{\mathrm{I}} &= \frac{2}{n+2} \cdot \frac{(2+\alpha)^n}{(1+\alpha)^n} \\
\frac{\mathrm{III}}{\mathrm{II}} &= \frac{3}{n+3} \cdot \frac{(3+\alpha)^n}{(2+\alpha)^n} \\
\frac{\mathrm{IV}}{\mathrm{III}} &= \frac{4}{n+4} \cdot \frac{(4+\alpha)^n}{(3+\alpha)^n} \\
\mathrm{etc.} &
\end{align*}
erhalten. Auf diese Weise wollen wir also die vorhergehenden in die folgenden einbeziehen, und weil $\mathrm{I}$
\[
	\Delta : n = \frac{1}{n+1}(1+\alpha)^n
\]
ist, werden wir für die Zahl $\mathrm{II}$ 
\[
	\Delta : n = \frac{1}{n+1}(1+\alpha)^n \cdot \frac{2}{n+2}\cdot \frac{(2+\alpha)^n}{(1+\alpha)^n}
\]
haben. Aus dieser geht weiter für die Zahl $\mathrm{III}$
\[
	\Delta : n = \frac{1}{n+1}(1+\alpha)^n \cdot \frac{2}{n+2} \cdot \frac{(2+\alpha)^n}{(1+\alpha)^n} \cdot \frac{3}{n+3} \cdot \frac{(3+\alpha)^n}{(2+\alpha)^n}
\]
hervor. Aus dieser wird auf die gleiche Weise für die Zahl $\mathrm{IV}$
\[
	\Delta : n = \frac{1}{n+1}(1+\alpha)^n \cdot \frac{2}{n+2} \cdot \frac{(2+\alpha)^n}{(1+\alpha)^n} \cdot \frac{3}{n+3} \cdot \frac{(3+\alpha)^n}{(2+\alpha)^n} \cdot \frac{4}{n+4} \cdot \frac{(4+\alpha)^n}{(3+\alpha)^n}
\]
hervorgehen.

\paragraph{§9}
Wenn wir also diese Ausdrücke ins Unendliche fortsetzen, werden wir schließlich die Wahrheit selbst erhalten; aber weil das Verhältnis des ersten Gliedes von den folgenden abweicht, wird sich das leicht zur Gleichmäßigkeit bringen lassen, indem  man den ganzen Ausdruck mit $\alpha^n$ erweitert; und so werden wir zum folgenden ins Unendliche laufende Produkt gelangen, wodurch der wahre Wert der Funktion $\Delta : n$ ausgedrückt werden wird; es wird nämlich
\[
	\Delta : n = \alpha^n \frac{1}{n+1}\ebinom{1+\alpha}{\alpha}^n \cdot \frac{2}{n+2}\ebinom{2+\alpha}{1+\alpha}^n \cdot \frac{3}{n+3}\ebinom{3+\alpha}{2+\alpha}^n \cdot \frac{4}{n+4}\ebinom{4+\alpha}{3+\alpha}^n
\]
\begin{flushleft}
\hspace{1em}$\mathrm{etc.}$
\end{flushleft}
sein, welcher Ausdruck daher natürlich einen bestimmten Wert angibt, welche Zahl auch immer, ob ganz oder gebrochen, für $n$ angenommen wird und deshalb sich diese Faktoren immer mehr der Einheit nähern; dass wird am besten aus der Form des infinitesimalen Faktors klar werden, welcher
\[
	\frac{i}{n+i}\cdot \ebinom{i+\alpha}{i-1+\alpha}^n
\]
ist, deren Wert wegen $i = \infty$ natürlich gleich $1$ ist, weil ja in Bezug auf $i$ die hinzugefügten $n$, $\alpha$, und $\alpha - 1$ verschwinden.

\paragraph{§10}
Diese Form stimmt schon überaus mit der, die wir eingangs erwähnt haben, überein, welche natürlich, wenn die Potenzen des Exponenten $n$ verbunden werden, auf diese Form
\[
	\Delta : n = \frac{1}{n+1}\ebinom{2}{1}^n \cdot \frac{2}{n+2}\ebinom{3}{2}^n \cdot \frac{3}{n+3}\ebinom{4}{3}^n \cdot \mathrm{etc.}
\]
zurückgeführt wird, auf welche die gerade gefundene übergeht, indem man $\alpha = 1$ setzt. Daraus wird eingesehen, dass die Formel, die wir nun gefunden haben, um vieles allgemeiner ist, während sich anstelle von $\alpha$ andere beliebige Zahlen annehmen lassen. Dannach besteht kein Zweifel, dass beide Formeln für alle Werte von $n$ selbst dieselben Werte ergeben. Das wurde zumindest aus den oben für die ganzen Zahlen ausgeführten Entwicklungen hinreichend klar gemacht, wo eines Beispiels wegen $\Delta : 3 = 1\cdot 2 \cdot 3$ hervorging, welchen Wert auch immer man für $\alpha$ angenommen hätte.

\paragraph{§11}
Dass aber die Größe des Buchstaben $\alpha$ natürlich nicht den Wert von $\Delta : n$ selbst betrifft, kann daraus leicht erkannt werden, dass alle Potenzen von $\alpha$ sich bis hin zur infinitesimalen Potenz sich gegenseitig aufheben; dann kann aber auch auf diese Weise gezeigt werden, wenn anstelle von $\alpha$ ein anderer beliebiger Wert $\beta$ geschrieben wird, dass dann in gleicher Weise 
\[
	\Delta : n = \beta^n\frac{1}{n+1}\ebinom{1+\beta}{\beta}^n \cdot \frac{2}{n+2}\ebinom{2+\beta}{1+\beta}^n \cdot \frac{3}{n+3}\ebinom{3+\beta}{2+\beta}^n \cdot \mathrm{etc.}
\]
sein wird. Es geht der Quotient, der aus der Division der einen Formel durch die andere entsteht,
\[
	1 = \ebinom{\alpha}{\beta}^n\left[\frac{(\alpha + 1)\beta}{\alpha(\beta+1)}\right] ^n \cdot \left[\frac{(\alpha + 2)(\beta + 1)}{(\alpha + 1)(\beta + 2)}\right] ^n \cdot \mathrm{etc.}
\]
hervor, woher, wenn die $n$te-Wurzel gezogen wird,
\[
	1 = \frac{\alpha}{\beta} \cdot \frac{(\alpha + 1)\beta}{\alpha(\beta + 1)} \cdot \frac{(\alpha + 2)(\beta + 1)}{(\alpha + 1)(\beta + 2)} \cdot \frac{(\alpha + 3)(\beta + 2)}{(\alpha + 2)(\beta + 3)} \cdot \mathrm{etc.}
\]
hervorgeht, welcher Bruch ins Unendliche fortgesetzt natürlich zu $\frac{\alpha + i}{\beta + i}$ reduziert wird, dessen Wert natürlich wegen des unendlichen $i$ die Einheit ist; und so wurde gezeigt, dass diese beiden Werte von $\Delta : n$ selbst einander gleich sind. Derselbe Beweis wird aber auch für die folgende Entwicklung gelten, wo wir alle hypergeometrischen Reihen im Allgemeinen betrachten wollen.
\section*{Erörterung der Analysis für die allgemeine hypergeometrische Reihe}
\paragraph{§13}
Wir wollen nun auf die gleiche Weise diese allgemeine hypergeometrische Reihe
\[
	a,~ a(a+b),~ a(a+b)(a+2b),~ a(a+b)(a+2b)(a+3b),~ \mathrm{etc.}
\]
betrachten, deren Term, der dem Index $n$ entspricht und natürlich aus $n$ Faktoren zusammengesetzt ist,
\[
	a(a+b)(a+2b)(a+3b) \cdots (a+(n-1)b) = \Delta : n
\]
gesetzt werde, weil er ja wegen der konstanten Buchstaben $a$ und $b$ wie eine Funktion der variablen Größe $n$ gesehen werden kann. Daher wird also aus der Gestalt der Aufbau
\begin{align*}
\Delta : (n+1) &= \Delta : n(a+bn) \\
\Delta : (n+2) &= \Delta : n(a+bn)(a+(n+1)b) \\
\Delta : (n+3) &= \Delta : n(a+bn)(a+(n+1)b)(a+(n+2)b) \\
\mathrm{etc.} &
\end{align*}
sein; und es wird sogar, wenn $i$ eine unendliche Zahl bezeichnet,
\[
	\Delta : (n+i) = \Delta : n(a+nb)(a+(n+1)b) \cdots (a+(n+i-1)b)
\]
sein.

\paragraph{§14}
Aus der Gestalt der Form selbst aber ist
\[
	\Delta : i = a(a+b)(a+2b) \cdots (a+(i-1)b),
\]
woher für die folgenden
\begin{align*}
\Delta : (i+1) &= \Delta : i(a+bi) \\
\Delta : (i+2) &= \Delta : i(a+bi)(a+(i+1)b) \\
\Delta : (i+3) &= \Delta : i(a+bi)(a+(i+1)b)(a+(i+2)b) \\
\mathrm{etc.} &
\end{align*}
werden wird, wo die über die neuen hinzukommenden Faktoren als untereinander gleich angesehen werden können, so dass
\[
	\Delta : (i+3) = \Delta : i(a+ib)^3
\]
ist, woher wir für eine unbestimmte Zahl $n$
\[
	\Delta : (i+n) = \Delta : i(a+ib)^n
\]
haben werden und sogar, weil ja mit gleichem Recht anstelle von $a+ib$ auch $a + b + ib$ oder $a+2b+ib$ hätte geschrieben werden können, können wir allgemein
\[
	\Delta : (i+n) = \Delta : i(\alpha + ib)^n
\]
schreiben, während $\alpha$ eine beliebige endliche Größe bezeichnet, die in Bezug auf $ib$ verschwindet. 

\paragraph{§15}
Weil wir daher ja einen doppelten Ausdruck für dieselbe Funktion $\Delta : (i+n)$ erhalten haben, werden wir, nachdem ein Vergleich aufgestellt wurde, daher 
\[
	\Delta : n = \frac{1 : i(\alpha + ib)^n}{(a+nb)(a+(n+1)b) \cdots (a+(n+i-1)b)}
\]
finden, wo wir im Nenner $i$ Faktoren haben werden. Wir wollen also anstelle von $\Delta : i$ wieder das Produkt selbst einsetzen, welches in gleicher Weise aus $i$ Faktoren besteht, woher der folgende Wert resultieren wird:
\[
	\Delta : n = \frac{a}{a+nb} \cdot \frac{a+b}{a+(n+1)b} \cdot \frac{a+2b}{a+(n+2)b} \cdots \frac{a+(i-1)b}{a+(n+i-1)b} (\alpha + ib)^n
\]
und dieser ist der wahre Wert von $\Delta : n$ selbst, solange für $i$ eine unendlich große Zahl angenommen wird, ohne Rücksicht auf den Index $n$, ob er eine ganze Zahl oder eine gebrochene oder sogar eine Wurzel ist.

\paragraph{§16}
Aus diesen sieht man ein, dass, wenn wir anstelle von $i$ endliche Werte annehmen, der Fehler umso kleiner wird, je größer die genommene Zahl $i$ war; damit diese Annäherung an den wahren Wert umso besser erkannt werden kann, wollen wir anstelle von $i$ der Reihe nach die Zahlen $ 1,\, 2,\, 3,\, 4,\, \mathrm{etc.} $ schreiben und die daher entstehenden Ausdrücke mit den Zeichen $\mathrm{I},\, \mathrm{II},\, \mathrm{III},\, \mathrm{etc.} $ bezeichnen und es wird
\begin{align*}
\mathrm{I.}~~~&\frac{a}{a+nb}(\alpha + b)^n \\
\mathrm{II.}~~&\frac{a}{a+nb} \cdot \frac{a+b}{a+(n+1)b} \cdot (\alpha + 2b)^n \\
\mathrm{III.}~&\frac{a}{a+nb} \cdot \frac{a+b}{a+(n+1)b} \cdot \frac{a+2b}{a+(n+2)b} \cdot (\alpha + 3b)^n \\
\mathrm{IV.}~~&\frac{a}{a+nb} \cdot \frac{a+b}{a+(n+1)b} \cdot \frac{a+2b}{a+(n+2)b} \cdot \frac{a+3b}{a+(n+3)b}(\alpha + 4b)^n \\
\mathrm{etc.} &
\end{align*}
hervorgehen. Diese berechnet man aber weiter zu
\begin{align*}
\frac{\mathrm{II}}{\mathrm{I}} &= \frac{a+b}{a+(n+1)b} \cdot \ebinom{\alpha + 2b}{\alpha + b}^n \\
\frac{\mathrm{III}}{\mathrm{II}} &= \frac{a+2b}{a+(n+2)b} \cdot \ebinom{\alpha + 3b}{\alpha + 2b}^n \\
\frac{\mathrm{IV}}{\mathrm{III}} &= \frac{a+3b}{a+(n+3)b} \cdot \ebinom{\alpha + 4b}{\alpha + 3b}^n \\
\mathrm{etc.} &
\end{align*}
Wenn wir also auf diese Weise ins Unendliche fortschreiten, werden wir zum wahren Wert von $\Delta : n $ selbst gelangen, der aus den folgenden miteinander zu multiplizierenden Faktoren bestehen wird, nachdem wir natürlich den ersten Faktor auf die Form der folgenden gebracht haben:
\begin{align*}
	\Delta : n =& \alpha^n\frac{a}{a+nb}\ebinom{\alpha + b}{\alpha}^n \cdot \frac{a+b}{a+(n+1)b}\ebinom{\alpha + 2b}{\alpha + b}^n \\
	 &\cdot\frac{a+2b}{a+(n+2)b}\ebinom{\alpha + 3b}{\alpha + 2b}^n \cdot \frac{a+3b}{a+(n+3)b}\ebinom{\alpha + 4b}{\alpha + 3b}^n \cdot \mathrm{etc.}
\end{align*}

\paragraph{§17}
Auf diese Weise haben wir also für $\Delta : n$ ein ins Unendliche laufende Produkt erhalten, dessen einzelne Terme nach einem hinreichend regelmäßigen Gesetz vorgehen; dort sollte besonders bemerkt werden, dass die einzelnen ganzen Faktoren oder die Glieder sich immer mehr der Wahrheit nähern, weil ja das infinitesimale Glied
\[
	\frac{a+(i-1)b}{a+(n+i-1)b}\ebinom{\alpha + ib}{\alpha + (i-1)b}^n
\]
sein wird, welcher Ausdruck, nachdem die Teile weggelassen wurden, die in Bezug auf einen unendlichen verschwinden, natürlich die Einheit wird. Dann haben wir aber schon beobachtet, dass die Größe $\alpha$ völlig unserem Belieben überlassen ist und daher der Wert $\Delta : n$ nicht verändert wird, woher in jedem Fall er sich so annehmen lassen wird, dass die Rechnung angenehmer wird; deshalb wird es gewiss der Mühen Wert sein, beobachtet zu haben, dass für Reihen solcher Art der allgemeine Term in einer um Vieles allgemeineren Form dargestellt werden kann als der, die wir wie eingangs angeführt haben, welche natürlich aus der gegenwärtigen dadurch entsteht, indem man $\alpha = a$ setzt.

\boldmath\section*{Anwendung dieser allgemeinen Formel beim Fall $n=\frac{1}{2}$}\unboldmath
\paragraph{§18} 
Man sieht leicht ein, dass der für $\Delta : n$ gefundene Ausdruck einen besonders großen Nutzen leisten kann, wann immer man Terme der Reihe wünscht, deren Indizes gebrochene Zahlen sind, weil ja die Terme, die den ganzen Zahlen entsprechen, per se bekannt sind. Wir wollen also zuerst den Term unserer Reihe suchen, der dem Index $n=\frac{1}{2}$ entspricht, welcher also durch $\Delta : \frac{1}{2}$ ausgedrückt werden wird, so dass
\[
	\Delta : \frac{1}{2} = \sqrt{\alpha}\frac{a}{a+\frac{1}{2}b}\sqrt{\frac{\alpha + b}{\alpha}} \cdot \frac{a+b}{a+\frac{3}{2}b}\sqrt{\frac{\alpha + 2b}{\alpha + b}}\cdot \frac{a+2b}{a+\frac{5}{2}b}\sqrt{\frac{\alpha + 3b}{\alpha + 2b}} ~~ \mathrm{etc.}
\]
bis ins Unendliche geht. Nachdem aber dieser Wert gefunden wurde, werden leicht alle Zwischenterme, die von zwei benachbarten gleichweit entfernt sind, bekannt; es wird nämlich
\[
	\Delta : 1\frac{1}{2} = \Delta : \frac{1}{2}(a + \frac{1}{2}b)
\]
sein, welcher Term zwischen den ersten $a$ und den zweiten $a(a+b)$ in die Mitte fällt; auf ähnliche Weise wird
\[
	\Delta : 2\frac{1}{2} = \Delta : \frac{1}{2}(a+\frac{1}{2}b)(a+\frac{3}{2}b)
\]
sein, welcher (Term) zwischen den zweiten und dritten genau in die Mitte fällt. Außerdem wird aber
\begin{align*}
\Delta : 3\frac{1}{2} &= \Delta : \frac{1}{2} (a+\frac{1}{2}b)(a + \frac{3}{2}b)(a+\frac{5}{2}b) \\
\Delta : 4\frac{1}{2} &= \Delta : \frac{1}{2} (a+\frac{1}{2}b)(a+\frac{3}{2}b)(a+\frac{5}{2}b)(a+\frac{7}{2}b) \\
\mathrm{etc.} &
\end{align*}
sein. Nun wird aber hier für gewöhnlich gefragt, wie diese ins Unendliche laufenden Produkte auf endlichen Ausdrücke zurückgeführt werden sollen, weil ja jene ins Unendliche erstreckten Produkte nur benutzt werden können, um den wahren Wert von $\Delta : \frac{1}{2}$ selbst näherungsweise zu finden. Besonders aber will man für gewöhnlich die Art der transzendenten Größen wissen, zu welcher dieser Wert $\Delta : \frac{1}{2}$ zu zählen ist; das kann aber durch das, was von mir vielerorts über unendliche Produkte dieser Art ausgelegt wurde, nicht geleistet werden. Vor allem ist es aber nötig, dass man die Wurzelfaktoren aus der Rechnung rauswirft, was durch Quadrieren passiert, woher wir
\[
	(\Delta : \frac{1}{2})^2 = \alpha \frac{aa}{(a+\frac{1}{2}b)^2} \cdot \frac{\alpha + b}{\alpha} \ebinom{a+b}{a+\frac{3}{2}b}^2 \frac{\alpha + 2b}{\alpha + b} \ebinom{a+2b}{a+\frac{5}{2}b}^2 \frac{\alpha + 3b}{\alpha + 2b} ~~ \mathrm{etc.}
\]
haben werden.

\paragraph{§19}
Weil aber hier der Buchstabe $\alpha$ von unserem Belieben abhängt, wollen wir ihn so annehmen, dass die Anzahl der Faktoren in den einzelnen Gliedern vermindert wird, was, indem man $\alpha = a$ nimmt, passieren wird; dann haben wir nämlich
\begin{align*}
	(\Delta : \frac{1}{2})^2 =& a\frac{a(a+b)}{(a+\frac{1}{2}b)(a+\frac{1}{2}b)} \cdot \frac{(a+b)(a+2b)}{(a+\frac{3}{2}b)(a+\frac{3}{2}b)} \\
	 &\cdot \frac{(a+2b)(a+3b)}{(a+\frac{5}{2}b)(a+\frac{5}{2}b)} \cdot \frac{(a+3b)(a+4b)}{(a+\frac{7}{2}b)(a+\frac{7}{2}b)} \cdot \mathrm{etc.}
\end{align*}
Damit wir aber daher die Partialbrüche aus den Nennern loswerden, wollen wird die einzelnen Faktoren des Zählers wie des Nenners verdoppeln, so dass diese Form
\[
	(\Delta : \frac{1}{2})^2 = a\frac{2a(2a+2b)}{(2a+b)(2a+b)} \cdot \frac{(2a+2b)(2a+4b)}{(2a+3b)(2a+3b)} \cdot \frac{(2a+4b)(2a+6b)}{(2a+5b)(2a+5b)} ~~ \mathrm{etc.}
\]
hervorgeht, wo die einzelnen Faktoren eines Gliedes für das folgende Glied einen Zuwachs von $2b$ erhalten. Nun kann aber diese Form leicht auf endliche Ausdrücke zurückgeführt werden durch das was überall (von mir) zerstreut erläutert worden ist.

\paragraph{§20}
Wenn nämlich mit den Buchstaben $P$ und $Q$ diese Integralformeln
\[
	P = \int\frac{x^{p-1}\partial x}{(1-x^n)^{1-\frac{m}{n}}} ~~ \mathrm{und} ~~ Q = \int \frac{x^{q-1}\partial x}{(1-x^n)^{1-\frac{m}{n}}}
\]
bezeichnet werden, welche Integrale natürlich so zu verstehen sind, dass sie von $x=0$ bis $x=1$ erstreckt werden, habe ich gezeigt, dass der Bruch $\frac{P}{Q}$ in das folgende unendliche Produkt verwandelt werden kann:
\[
	\frac{P}{Q} = \frac{q(m+p)}{p(m+q)}\cdot \frac{(q+n)(m+p+n)}{(p+n)(m+q+n)} \cdot \frac{(q+2n)(m+p+2n)}{(p+2n)(m+q+2n)}\cdot \mathrm{etc.},
\]
wo die einzelnen Faktoren eines jeden Gliedes immer um die Größe $n$ wachsen; daher ist sofort klar, damit diese Form auf diese, die uns vorgelegt ist, übergeht, $n=2b$ genommen werden muss; dann aber genügt es, die ersten Glieder auf beiden Seiten einander gleichzusetzen, was natürlich
\[
	\frac{q(m+p)}{p(m+q)} = \frac{2a(2a+2b)}{(2a+b)(2a+b)}
\]
ergibt; das kann aber gemacht werden, indem man $q=2a$ und $p=2a+b$ nimmt, dann aber $m=b$, nach dem Einsetzen welcher Werte der Bruch $\frac{P}{Q}$ mit $a$ multipliziert den Wert $(\Delta : \frac{1}{2})^2$ selbst, den wir suchen, auf endliche Weise ausdrückt.

\paragraph{§21}
Nachdem aber die gerade gefundene Substitution gemacht wurde, wird wegen $\frac{m}{n} = \frac{1}{2}$
\[
	P = \int\frac{x^{2a+b-1}\partial x}{\sqrt{1-x^{2b}}} ~~ \mathrm{und} ~~ Q = \int\frac{x^{2a-1}\partial x}{\sqrt{1-x^{2b}}}
\]
sein, welche Integrale natürlich immer von $x=0$ bis $x=1$ zu erstrecken sind; dadurch wird
\[
	(\Delta : \frac{1}{2})^2 = a\frac{P}{Q}
\]
sein, und daher wird, indem man die Wurzel zieht
\[
	\Delta : \frac{1}{2} = \sqrt{a\int\frac{x^{2a+b-1}\partial x}{\sqrt{1-x^{2b}}} : \int\frac{x^{2a-1}\partial x}{\sqrt{1-x^{2b}}}}
\]
sein, aus welcher Formel in einem Fall sofort klar werden wird, von welcher Art transzendenter Größen der gesuchte Wert $\Delta : \frac{1}{2}$ abhängt; es wird der Mühen Wert sein, das an einigen Beispielen zu illustrieren.

\section*{Beispiel 1}
\paragraph{§22}
Man setze $a=1$ und $b=1$, sodass die Wallis'sche hypergeometrische Reihe
\[
	1,~ 1 \cdot 2,~ 1 \cdot 2 \cdot 3,~ 1 \cdot 2 \cdot 3 \cdot 4,~ 1 \cdot 2 \cdot 3 \cdot 4 \cdot 5,~ \mathrm{etc.}
\]
selbst hervorgeht, deren Term, der dem Index $\frac{1}{2}$ entspricht und durch $\Delta : \frac{1}{2}$ bezeichnet wurde, gesucht wird. Durch die gefundene Formel wird also
\[
	\Delta : \frac{1}{2} = \sqrt{\int\frac{xx\partial x}{\sqrt{1-xx}} : \int\frac{x\partial x}{\sqrt{1-xx}}}
\]
sein. Es ist aber bekannt, nachdem diese Integrale von $x=0$ bis $x=1$ erstreckt wurden, dass zuerst
\[
	\int\frac{x\partial x}{\sqrt{1-xx}} = 1
\]
ist, dann aber
\[
	\int\frac{xx\partial x}{\sqrt{1-xx}} = \frac{1}{2}\int\frac{\partial x}{\sqrt{1-xx}} = \frac{\pi}{4}
\]
woher klar ist, dass
\[
	\Delta : \frac{1}{2} = \sqrt{\frac{\pi}{4}} = \frac{1}{2}\sqrt{\pi}
\]
sein wird; die übrigen Zwischenterme dieser Reihe werden aber
\begin{align*}
	\Delta : 1\frac{1}{2} &= \frac{1}{2}\cdot \frac{3}{2}\sqrt{\pi} \\
	\Delta : 2\frac{1}{2} &= \frac{1}{2}\cdot \frac{3}{2}\cdot \frac{5}{2}\sqrt{\pi} \\
	\Delta : 3\frac{1}{2} &= \frac{1}{2}\cdot \frac{3}{2}\cdot \frac{5}{2} \cdot \frac{7}{2}\sqrt{\pi}\\
	\mathrm{etc.} &
\end{align*}
sein, woher klar ist, dass der in der Reihe vorhandene Term, der dem Index $-\frac{1}{2}$ entspricht,
\[
	\Delta : -\frac{1}{2} = \sqrt{\pi}
\]
sein wird, genauso wie von Wallis bemerkt worden ist.
\section*{Beispiel 2}
Man nehme $a=1$ und $b=2$, woher diese hypergeometrische Progression entsteht:
\[
	1,~ 1 \cdot 3,~ 1 \cdot 3 \cdot 5,~ 1 \cdot 3 \cdot 5 \cdot 7,~ 1 \cdot 3 \cdot 5 \cdot 7 \cdot 9,~ \mathrm{etc.} 
\]
Der andere Term, der dem Index $\frac{1}{2}$ entspricht, wird gesucht; er wird also
\[
	\Delta : \frac{1}{2} = \sqrt{\int\frac{x^3\partial x}{\sqrt{1-x^4}} : \int\frac{x \partial x}{\sqrt{1-x^4}}}
\]
sein. Wenn wir also hier anstelle von $xx$ gleich $y$ schreiben, werden wir
\[
	\int\frac{x^3 \partial x}{\sqrt{1-x^4}} = \frac{1}{2}\int \frac{y\partial y}{\sqrt{1-yy}}
\]
haben, deren Wert von $y=0$ bis $y=1$ erstreckt gleich $\frac{1}{2}$ ist; die andere Formel
\[
	\int \frac{x\partial x}{\sqrt{1-x^4}}
\]
aber geht über in 
\[
	\frac{1}{2}\int\frac{\partial y}{\sqrt{1-yy}} = \frac{1}{2} \cdot \frac{\pi}{2}.
\]
Nach Einsetzen dieser Werte wird
\[
	\Delta : \frac{1}{2} = \sqrt{\frac{2}{\pi}}
\]
sein, welcher Wert auch von der Quadratur des Kreises abhängt. Dann aber werden die Zwischenterme die folgenden sein:
\begin{align*}
\Delta : 1\frac{1}{2} &= 2\sqrt{\frac{2}{\pi}} \\
\Delta : 2\frac{1}{2} &= 2\cdot 4\sqrt{\frac{2}{\pi}} \\
\Delta : 3\frac{1}{2} &= 2\cdot 4\cdot 6 \sqrt{\frac{2}{\pi}} \\
\Delta : 4\frac{1}{2} &= 2\cdot 4\cdot 6 \cdot 8\sqrt{\frac{2}{\pi}} \\
\mathrm{etc.} &
\end{align*}
und daher ist klar, dass der Term, der dem Index $-\frac{1}{2}$ entspricht, unendlich sein wird.

\boldmath\section*{Anwendung, um den Term dieser Reihen zu finden, dessen Index gleich $\frac{1}{3}$ ist}\unboldmath
\paragraph{§24}
Wir wollen also hier $n=\frac{1}{3}$ setzen und unsere gefundene, allgemeine Formel wird uns
\[
	\Delta : \frac{1}{3} = \sqrt[3]{\alpha} \cdot \frac{a}{a+\frac{1}{3}b}\sqrt[3]{\frac{\alpha + b}{\alpha}} \cdot \frac{a+b}{a+\frac{4}{3}b}\sqrt[3]{\frac{\alpha + 2b}{\alpha + b}} \cdot \frac{a+2b}{a+\frac{7}{3}b}\sqrt[3]{\frac{\alpha + 3b}{\alpha + 2b}} \cdot \mathrm{etc.}
\]
liefern. Daher wird also durch Kubieren
\[
	(\Delta : \frac{1}{3})^3 = \alpha \cdot \frac{a^3}{(a+\frac{1}{3}b)^3} \cdot \frac{\alpha + b}{\alpha} \cdot \frac{(a+b)^3}{(a+\frac{4}{3}b)^3} \cdot \frac{\alpha + 2b}{\alpha + b}  \cdot \mathrm{etc.}
\]
sein. Um gleich die Brüche loszuwerden, wollen wir $b=3c$ setzen und es wird
\[
	(\Delta : \frac{1}{3})^3 = \alpha \cdot \frac{a^3}{(a+c)^3} \cdot \frac{\alpha + 3c}{\alpha} \cdot \frac{(a+3c)^3}{(a+4c)^3} \cdot \frac{\alpha + 6c}{\alpha + 3c} \cdot \frac{(a+6c)^3}{(a+7c)^3} \cdot \frac{\alpha + 9c}{\alpha + 6c} \cdot \mathrm{etc.}
\]
werden, wo der Zuwachs $3c$ ist, welchen die einzelnen Faktoren haben, während wir von einem Glied aus zum folgenden gehen. Wenn wir also $\alpha = 0$ setzen, werden wird zur folgenden einfacheren Form gelangen:
\begin{align*}
	(\Delta : \frac{1}{3})^3 =& a \cdot \frac{aa(a+3c)}{(a+c)(a+c)(a+c)} \cdot \frac{(a+3c)(a+3c)(a+6c)}{(a+4c)(a+4c)(a+4c)} \\
	&\cdot \frac{(a+6c)(a+6c)(a+9c)}{(a+7c)(a+7c)(a+7c)} \cdot \mathrm{etc.}
\end{align*}

\paragraph{§25}
Weil ja hier in jedem Glied drei Faktoren auftauchen, lässt sich ein Vergleich mit der für $\frac{P}{Q}$ beschaffenen From nicht sofort aufstellen. Aber hier muss man zwei Brüche $\frac{P}{Q}$ und $\frac{P'}{Q'}$ solcher Art, deren Produkt man der gefundenen Form gleichsetze, zur Hilfe nehmen; und weil ja unser erstes Glied gleich
\[
	\frac{aa(a+3c)}{(a+c)(a+c)(a+c)}
\]
ist, erzeugen die zwei ersten Glieder, die auch bei jener Multiplikation entstehen, vier Faktoren; wir wollen bei den einzelnen Gliedern stets mit dem neuen Faktor $f$ erweitern, sodass sie in $2$ Teile gespalten werden können, welche für das erste Glied
\[
	\frac{aa}{(a+c)f}\cdot \frac{f(a+3c)}{(a+c)(a+c)}
\]
seien, und nun wollen wir beide Teile mit
\[
	\frac{q(m+p)}{p(m+q)}
\]
vergleichen. Für den ersten Teil aber wollen wir
\[
	q = a ~~ \mathrm{und} ~~ p = f
\]
setzen und es wird
\[
	m + p = m + f = a
\]
werden und
\[
	m + q = m + a = a + c,
\]
woher man $m=c$ und $f=a-c$ berechnet; dann aber, indem man zu folgenden Gliedern vorgeht, wird $n=3c$ werden.

\paragraph{§26}
Wenn wir also gleich die einzelnen Glieder unseres Ausdrucks in zwei solche Teile auflösen, wollen wir, nachdem der neue Buchstabe $f=a-c$, der in gleicher Weise bei den folgenden Gliedern einen Zuwachs von $3c$ erhalten wird, alle ersten Teile getrennt betrachten, deren Produkt dem Bruch $\frac{P}{Q}$ gleich werden wird; und es wird aus den schon gefundenen Werten
\[
	P = \int\frac{x^{a-c-1}\partial x}{\sqrt[3]{(1-x^{3c})^2}} ~~ \mathrm{und} ~~ Q = \int\frac{x^{a-1}\partial x}{\sqrt[3]{(1-x^{3c})^2}}
\]
sein.

\paragraph{§27}
Für die weiteren Teile aber wollen wir den Bruch $\frac{P'}{Q'}$ zur Hilfe nehmen, während wir auch die Kleinbuchstaben $p$, $q$ und $m$ mit dem gleichen Zeichen versehen wollen. Daher wird also der Vergleich unserer ersten Glieder
\[
	\frac{q'(m'+p')}{p'(m'+q')} = \frac{f(a+3c)}{(a+c)(a+c)}
\]
ergeben; deswegen wollen wir
\[
	q' = f = a - c ~~ \mathrm{und} ~~ p' = a + c
\]
nehmen; dann wird aber
\[
	m' + p' = m' + a + c = a + 3c
\]
sein und
\[
	m' + q' = m' + a - c = a + c,
\]
wo auf jeder der beiden Seiten $m'=2c$ wird; dann aber bleibt wie zuvor $n=3c$, woraus diese Formeln so bestimmt werden:
\[
	P'=\int\frac{x^{a+c-1}\partial x}{\sqrt[3]{1-x^{3c}}} ~~ \mathrm{und} ~~ Q' = \int\frac{x^{a-c-1}\partial x}{\sqrt[3]{1-x^{3c}}}
\]

\paragraph{§28}
Weil also der Bruch $\frac{P}{Q}$ ein Produkt aller ersten Teile ausdrückt, aber in der Tat $\frac{P'}{Q'}$ ein Produkt aller anderen Teile, werden wir
\[
	(\Delta : \frac{1}{3})^3 = a\frac{P}{Q}\cdot \frac{P'}{Q'}
\]
haben und daher
\[
	(\Delta : \frac{1}{3}) = \sqrt[3]{\frac{aPP'}{QQ'}},
\]
und so wird man, um diesen interpolierten Term $\Delta : \frac{1}{3}$ zu bestimmen $4$ Integralformeln brauchen, zwischen denen man im Allgemeinen keine Relation erkennt; es wird nämlich, nachdem die Substitution gemacht wurde,
\[
\Delta : \frac{1}{3} = \sqrt[3]{a\int\frac{x^{a-c-1}\partial x}{\sqrt[3]{(1-x^{3c})^2}} \cdot \int\frac{x^{a+c-1}\partial x}{\sqrt[3]{(1-x^{3c})}}} : \sqrt[3]{\int\frac{x^{a-1}\partial x}{\sqrt[3]{(1-x^{3c})^2}} \cdot \int \frac{x^{a-c-1}\partial x}{\sqrt[3]{(q-x^{3c})}}}
\]
sein, wo man sich erinnern sollte, dass anstelle des Buchstabens $b$ hier $3c$ geschrieben wurde, sodass $c=\frac{1}{3}b$ ist. Es soll genügen, diesen Ausdruck am Beispiel der Wallis'schen Reihe illustriert zu haben.

\section*{Beispiel}
\paragraph{§29}
Es sei also $a=1$ und $b=1$ und daher $c=\frac{1}{3}$; und daher werden wir die vier Integralformeln
\[
	P = \int\frac{x^{-\frac{1}{3}}\partial x}{\sqrt[3]{(1-x)^2}} ~~ \mathrm{und} ~~ Q = \int\frac{\partial x}{\sqrt[3]{(1-x)^2}}
\]
\[
	P' = \int\frac{x^{\frac{1}{3}}\partial x}{\sqrt[3]{1-x}} ~~ \mathrm{und} ~~ Q' = \int\frac{x^{-\frac{1}{3}}\partial x}{\sqrt[3]{1-x}}
\]
sein; um diese Formeln von gebrochenen Exponenten zu befreien, setze man $x=y^3$ und es wird daher
\[
	P = 3\int\frac{y\partial y}{\sqrt[3]{(1-y^3)^2}} ~~ \mathrm{und} ~~ Q = 3\int\frac{yy\partial y}{\sqrt[3]{(1-y^3)^2}}
\]
\[
	P' = 3\int\frac{y^3 \partial y}{\sqrt[3]{1-y^3}} ~~ \mathrm{und} ~~ Q' = 3\int\frac{y\partial y}{\sqrt[3]{1-y^3}}
\]
sein, aus welchen Werten
\[
	\Delta : \frac{1}{3} = \sqrt[3]{\frac{aPP'}{QQ'}}
\]
wird, wo man bemerke, dass
\[
	\int\frac{y^3\partial y}{\sqrt[3]{1-y^3}} = \frac{1}{3}\int\frac{\partial y}{\sqrt[3]{1-y^3}}
\]
ist.

\paragraph{§30}
Wir wollen gleich diese einzelnen Integralformeln sorgfältiger entwickeln; und zwar kann zuerst die Formel
\[
	\int\frac{\partial y}{\sqrt[3]{1-y^3}}
\]
auf den Kreis zurückgeführt werden. Nachdem nämlich
\[
	\frac{y}{\sqrt[3]{1-y^3}} = z
\]
gesetzt wurde, wird unsere Formel $\frac{z\partial y}{y}$ sein; dann wird aber
\[
	y^3 = \frac{z^3}{1+z^3} ~~ \mathrm{und} ~~ 3\log{y} = 3\log{z} - \log{(1+z^3)}
\]
sein und daher
\[
	\frac{\partial y}{y} = \frac{\partial z}{z} - \frac{zz\partial z}{1+z^3} = \frac{\partial z}{z(1+z^3)}
\]
und so wird unsere Formel gleich
\[
	\int \frac{\partial z}{1+z^3}
\]
werden, deren Integral von $y=0$ bis $y=1$, dass heißt von $z=0$ bis $z=\infty$ gleich
\[
	\frac{\pi}{3\sin{\frac{\pi}{3}}} = \frac{2\pi}{3\sqrt{3}}
\]
ist, und so wird $P' = \frac{2\pi}{3\sqrt{3}}$ sein. Aber für $P$ wird die Formel
\[
	\int \frac{y\partial y}{\sqrt[3]{(1-y^3)^2}}
\]
durch dieselbe Substitution
\[
	\frac{y}{\sqrt[3]{1-y^3}} = z
\]
auf diese
\[
	\int \frac{zz\partial y}{y} = \int \frac{z\partial z}{1+z^3},
\]
zurückgeführt, deren Integral
\[
	\frac{\pi}{3\sin{\frac{2\pi}{3}}} = \frac{2\pi}{3\sqrt{3}}
\]
ist, woher $P = \frac{2\pi}{\sqrt{3}}$ wird. Weiter aber ist für $Q$
\[
	\int \frac{yy\partial y}{\sqrt[3]{(1-y^3)^2}} = 1 - \sqrt[3]{1-y^3},
\]
woher für $y=1$ gesetzt $Q=3$ werden wird. Schließlich lässt sich die Formel
\[
	Q' = 3\int\frac{y\partial y}{\sqrt[3]{1-y^3}}
\]
auf keine Weise auf bekannte Größen zurückführen, sondern involviert eine einzige Quadratur. Aus diesen Formeln bildet man dann
\[
	\Delta : \frac{1}{3} = \sqrt[3]{\frac{4\pi\pi}{81\int\frac{y\partial y}{\sqrt[3]{1-y^3}}}},
\]
aus welchem Wert man weiter die folgenden schließt:
\begin{align*}
\Delta : 1\frac{1}{3} &= \frac{4}{3}\Delta : \frac{1}{3} \\
\Delta : 2\frac{1}{3} &= \frac{4}{3} \cdot \frac{7}{3} \Delta : \frac{1}{3} \\
\Delta : 3\frac{1}{3} &= \frac{4}{3} \cdot \frac{7}{3} \cdot \frac{10}{3} \Delta : \frac{1}{3} \\
\mathrm{etc.} &
\end{align*}
Daher sieht man aber leicht ein wie die Untersuchung gehandhabt werden muss, wenn anstelle von $n$ andere Brüche gesetzt werden.
\section*{Zusammenfassung und Schlussbemerkung}
\paragraph{§31}
So wie wir hier für die Reihe
\[
	a,~ a(a+b),~ a(a+b)(a+2b),~ a(a+b)(a+2b)(a+3b),~ \mathrm{etc.}
\]
den Term $\Delta : n$, der dem unbestimmten Index $n$ entspricht, so ausgedrückt gefunden haben, dass
\[
	\Delta : n = \alpha^n \cdot \frac{a}{a+nb}\ebinom{\alpha + b}{\alpha}^n \cdot \frac{a+b}{a+(n+1)b}\ebinom{\alpha + 2b}{\alpha + b}^n ~~ \mathrm{etc.}
\]
ist, wenn wir anstelle von $a$ eine andere beliebige Zahl $c$ annehmen, dass die Reihe
\[
c,~ c(c+b),~ c(c+b)(c+2b),~ c(c+b)(c+2b)(c+3b),~ \mathrm{etc.}
\]
ist und wir ihren Term, der dem Index $n$ entspricht mit $\Gamma : n$ bezeichnen, dann wird auf die gleiche Weise
\[
	\Gamma : n = \alpha^n\cdot\frac{c}{c+nb}\ebinom{\alpha + b}{\alpha}^n \cdot \frac{c+b}{c+(n+1)b}\ebinom{\alpha + 2b}{\alpha + b}^n ~~ \mathrm{etc.}
\]
sein, wo freilich $\alpha$ einen anderen beliebigen Wert bezeichnen könnte als in der ersten  Formel. Wenn wir also nun die letzte Reihe durch die erste teilen, wird daher die folgende Reihe entstehen:
\[
\frac{a}{c},~ \frac{a(a+b)}{c(c+b)},~ \frac{a(a+b)(a+2b)}{c(c+b)(c+2b)},~ \frac{a(a+b)(a+2b)(a+3b)}{c(c+b)(c+2b)(c+3b)},~ \mathrm{etc.}
\]
und es ist klar, dass der dem Index $n$ entsprechende Term $\frac{\Delta : n}{\Gamma : n}$ sein wird; wenn wir daher bei beiden für $\alpha$ dieselbe Zahl annehmen, werden alle Potenzen des Exponenten $n$ sich gegenseitig aufheben, sodass für diese Reihe der allgemeine oder dem Index $n$ entsprechende Term gleich
\[
\frac{a(c+nb)}{(a+nb)c} \cdot \frac{(a+b)(c+(n+1)b)}{(a+(n+1)b)(c+b)} \cdot \frac{(a+2b)(c+(n+2)b)}{(a+(n+2)b)(c+2b)} ~~ \mathrm{etc.}
\]
ist, wo also die Interpolation ohne jede Schwierigkeit aufgestellt werden kann. Ja, es kann sogar im Allgemeinen dieser allgemeine Term selbst angenehm durch den Bruch $\frac{P}{Q}$ ausgedrückt werden, indem man $q=a$, $p=c$, $m=nb$ nimmt, sodass der sukzessive Zuwachs, welcher $n$ war, nun $b$ ist, woher sich die beiden Integralformeln so verhalten werden
\[
	P = \int \frac{x^{c-1}\partial x}{(1-x^b)^{1-n}} ~~ \mathrm{und} ~~Q = \int\frac{x^{a-1}\partial x}{(1-x^b)^{1-n}};
\]
und so wird sich in diesen Fällen der allgemeine Term immer durch zwei Integralformeln und daher auf endliche und bestimmte Weise ausdrücken lassen und die Interpolation erfordert keine neuen Quadraturen, so wie es in den oben betrachteten Beispielen der Fall war.
\end{document}